\documentclass[11pt,reqno]{amsart}
\usepackage{latexsym,amscd,amssymb}
\usepackage{times,mathptm}
\newtheorem{thm}{Theorem}
\newtheorem*{thm-nn}{Theorem}
\newtheorem*{conj}{Conjecture}
\newtheorem{lem}[thm]{Lemma}
\newtheorem{pro}[thm]{Proposition}
\newtheorem*{pro-nn}{Proposition}

\newtheorem*{cor-nn}{Corollary}
\newtheorem*{conj-nn}{Conjecture}

\theoremstyle{definition}

\theoremstyle{remark}

\newtheorem*{notation}{Notation}
\newtheorem{rem}{Remark}
\newtheorem*{ack}{Acknowledgement}

\newcommand{\calS}{\mathcal{S}}

\newcommand{\Q}{\mathbb{Q}}
\newcommand{\R}{\mathbb{R}}
\newcommand{\Z}{\mathbb{Z}}

\newcommand{\funcS}{\mathbb{S}}

\newcommand{\diff}{\mathrm{Diff}_+{\Sigma_g}}

\newcommand{\lgr}{\langle\gamma\rangle}

\begin{document}
\title[Integral Riemann-Roch formulae]{
Integral Riemann-Roch formulae \\ for cyclic subgroups of
mapping class groups}
\author[T. Akita]{Toshiyuki Akita}
\address{Department of Mathematics, Hokkaido University,
Sapporo, 060-0810 Japan}
\email{akita@math.sci.hokudai.ac.jp}
\author[N. Kawazumi]{Nariya Kawazumi}
\address{Department of Mathematical Sciences,
The University of Tokyo, %3-8-1 Komaba, Meguro, 
Tokyo, 153-8914 Japan}
\email{kawazumi@ms.u-tokyo.ac.jp}
%\address{Fax: 011-717-9303}

%\thanks{This research was partially supported by the Ministry
%of Education, Science, Sports and Culture, Grant-in-Aid for  
%Scientific Research (C) (No. 17540056).}

%Grand-in-Aid for Encouragement of
%Young Scientists (No. 12740030), the Ministry of Education, Science,
%Sports and Culture.}
%\addtolength{\baselineskip}{0.4mm}

\setlength{\topmargin}{-15pt}
\setlength{\oddsidemargin}{25pt}
\setlength{\evensidemargin}{25pt}
\setlength{\textwidth}{400pt}
\setlength{\textheight}{640pt}

\addtolength{\parskip}{2mm}
\setlength{\parindent}{0.0mm}

\maketitle
\begin{abstract}
The first author conjectured % in \cite{nilpotency}
certain relations for
Morita-Mumford classes and Newton classes in the
integral cohomology of mapping class groups
(integral Riemann-Roch formulae).
In this paper, the conjecture is verified for cyclic
subgroups of mapping class groups.
\end{abstract}
\section{Introduction}
Let $\Sigma_g$ be a closed oriented surface of genus
$g$ and $\diff$ the group of orientation preserving
diffeomorphisms of $\Sigma_g$ equipped with the
$C^{\infty}$-topology.
The mapping class group $\Gamma_g$ is the group of
path components of $\diff$.
Throughout this paper we assume $g\geq 2$.
There are two important families of integral
cohomology classes of $\Gamma_g$, Morita-Mumford classes
$e_k$ and Chern classes
$c_k$.
Newton classes $s_k$ are integral polynomials of $c_i$'s
(see Section \ref{sec-def} for the precise definitions).
Over the rationals, these classes are related by the
Grothendieck-Riemann-Roch theorem (or the Atiyah-Singer
index theorem) which implies that
\[
\frac{B_{2k}}{2k}e_{2k-1}=s_{2k-1}\in H^{4k-2}(\Gamma_g;\Q)
\]
for all $k\geq 1$, where $B_{2k}$ is the $2k$-th Bernoulli number
(see
\cite{MR89e:57022,MR85j:14046}).
The first author conjectured in \cite{nilpotency} that
similar relations hold over the {\em integers},
as stated in the below.
\begin{conj}[Integral Riemann-Roch formulae for $\Gamma_g$]
%Let $N_{2k}'$ and $D_{2k}'$ be the numerator and 
%the denominator of $B_{2k}/2k$, respectively. Then
Let $N_{2k}',D_{2k}'$ be coprime integers satisfying 
$B_{2k}/2k=N_{2k}'/D_{2k}'$. Then
%\[
$N_{2k}'e_{2k-1}=D_{2k}'s_{2k-1}\in H^{4k-2}(\Gamma_g;\Z)$
%\]
holds for all $g\geq 2$ and $k\geq 1$.
\end{conj}
The conjecture is known to be affirmative for $k=1$ 
(i.e. $e_1=12s_1
\in H^2(\Gamma_g;\Z)$ for $g\geq 2$), 
since $H^2(\Gamma_g;\Z)\cong\Z$
for $g\geq 3$ as was proved by Harer
\cite{harer-2nd} (see \cite{nilpotency} for the
case $g=2$).
Moreover, the second author
showed that a slightly weaker version of the
conjecture holds for hyperelliptic mapping class groups
\cite{hyperell}.
The purpose of this paper is to verify the conjecture for
cyclic subgroups of $\Gamma_g$:
\begin{thm}\label{thm-main}
Let $\lgr$ be a cyclic subgroup of $\Gamma_g$ generated
by $\gamma\in\Gamma_g$. Then
\[
N_{2k}'e_{2k-1}(\gamma)=D_{2k}'s_{2k-1}(\gamma)
\in H^{4k-2}(\lgr;\Z),
\]
where $e_{2k-1}(\gamma),s_{2k-1}(\gamma)$ are the restriction
of $e_{2k-1},s_{2k-1}$ to $\lgr$.
\end{thm}
The rest of the paper is organized as follows.
In Section \ref{sec-def}
we will recall the definitions of
Morita-Mumford classes $e_k$ and Newton classes $s_k$
of $\Gamma_g$.
In Section \ref{sec-fpd} we will introduce the fixed
point data of an element $\gamma\in\Gamma_g$ of finite order,
and express $e_k(\gamma)$ and $s_k(\gamma)$ in terms of
the fixed point data of $\gamma$.
The expression of $e_k(\gamma)$ was obtained in \cite{AKU,MR1664755}.
The proof of Theorem \ref{thm-main}
will be given in Section \ref{sec-proof}. The key ingredient
for the proof of Theorem \ref{thm-main}
is a certain generalization of Voronoi's congruence
on Bernoulli numbers (Theorem \ref{thm-voronoi}),
which will be proved in
Section \ref{sec-Voronoi} by using a result of Porubsk\'y
\cite{porubsky,porubsky2}.

\begin{notation}
There are some conflicting notations of Bernoulli numbers.
We define $B_{2k}$ $(k\geq 1)$ by a power series expansion
\[
\frac{z}{e^z-1}+\frac{z}{2}=1+\sum_{k=1}^{\infty}
\frac{B_{2k}}{(2k)!}z^{2k}.
\]
Thus $B_2=1/6, B_4=-1/30,
B_6=1/42, B_8=-1/30, B_{10}=5/66$, and so on.
Our notation is consistent with \cite{IR}
but differs from
\cite{nilpotency,hyperell,MR89e:57022}.
\end{notation}

\section{Definitions}\label{sec-def}
In this section, we recall definitions of Morita-Mumford classes
and Newton classes of $\Gamma_g$.
Let $\pi:E\rightarrow B$ be an oriented $\Sigma_g$-bundle,
$T_{E/B}$ the tangent bundle along the fiber of $\pi$,
and $e\in H^2(E;\Z)$ the Euler class of $T_{E/B}$.
Then the $k$-th Morita-Mumford class $e_k\in H^{2k}(B;\Z)$
of $\pi$ is defined by
\[
e_k:=\pi_!(e^{k+1})
\]
where $\pi_!:H^*(E;\Z)\rightarrow H^{*-2}(B;\Z)$ is the Gysin
homomorphism (or the integration along the fiber).
The structure group of oriented $\Sigma_g$-bundles is $\diff$.
Hence passing to the universal $\Sigma_g$-bundle
\[
E\diff\times_{\diff}\Sigma_g\rightarrow B\diff
\]
we obtain the cohomology classes $e_k\in H^*(B\diff;\Z)$.
Now a result of Earle and Eells \cite{ee}
implies that the classifying space $B\diff$ is
homotopy equivalent to the Eilenberg-MacLane space
$B\Gamma_g=K(\Gamma_g,1)$ so that we obtain the
{\em $k$-th Morita-Mumford class
 $e_k\in H^{2k}(\Gamma_g;\Z)$
of} $\Gamma_g$.

Let $U(g)$ be the unitary group.
The $k$-th Newton class $s_k\in H^*(BU(g);\Z)$
is the characteristic class associated to the formal sum
$\sum_i x_i^k$.
Note that $s_1=c_1$ and $s_k=N_k(c_1,c_2,\dots,c_k)$
where $c_i\in H^*(BU(g);\Z)$ is the $i$-th Chern class
and $N_k$ is the $k$-th Newton polynomial.
The Newton classes are also called the integral Chern
character in the literature,
for $g+\sum_{k\geq 1}s_k/k!\in H^*(BU(g);\Q)$
coincides with the Chern character.

Now the natural action of $\Gamma_g$
on the first real homology
$H_1(\Sigma_g;\R)$ induces a homomorphism $\Gamma_g\rightarrow
Sp(2g,\R)$. The homomorphism yields a continuous map
\[
\eta:B\Gamma_g\rightarrow BU(g)
\]
of classifying spaces,
for the maximal compact subgroup of $Sp(2g,\R)$ is 
isomorphic to $U(g)$.
The {\em $k$-th Newton class $s_k\in H^{2k}(\Gamma_g;\Z)$ 
of} $\Gamma_g$ is defined to be
the pull-back of $s_k\in H^*(BU(g);\Z)$ by $\eta$
(we use the same symbol).
The {\em $k$-th Chern class $c_k\in H^{2k}(\Gamma_g;\Z)$ 
of} $\Gamma_g$ is defined similarly.

\begin{rem}
In \cite{nilpotency} and \cite{MR89e:57022},
Chern classes and Newton
classes of $\Gamma_g$ were defined in terms of the action of
$\Gamma_g$ on the cohomology $H^1(\Sigma_g;\R)$ rather than
the homology. Hence the $k$-th Newton class $s_k$ considered
in those papers is $(-1)^k$ times of ours.
\end{rem}
\begin{rem}\label{rem:holo}
Any oriented $\Sigma_g$-bundle $\pi:E\rightarrow B$
can be regarded as a continuous family of compact Riemann
surface of genus $g$.
Let $E_b$ be the fiber of of $b\in B$ and
$H^0(E_b,\Omega_{E_b}^1)$
the space of holomorphic 1-forms on $E_b$.
We obtain a $g$-dimensional complex vector bundle
\[
\xi:\bigcup_{b\in B} H^0(E_b,\Omega_{E_b}^1)\rightarrow B
\]
over $B$, which is called the Hodge bundle associated to $\pi$.
The Chern classes $c_k(\xi)$ and Newton classes $s_k(\xi)$
coincide with the pull-back of Chern classes and Newton
classes of $\Gamma_g$ by the
classifying map $B\rightarrow B\Gamma_g$ of $\pi$.
Indeed, the classifying map $B\rightarrow BU(g)$ of $\xi$
can be identified with the composition of
$B\rightarrow B\Gamma_g$
and $\eta:B\Gamma_g\rightarrow BU(g)$.
See \cite[p.~366]{hyperell} and
\cite[p.~553]{MR89e:57022} for precise.
\end{rem}

For later use, we introduce characteristic classes of
unitary representations of finite groups.
See \cite{thomas} for precise.
Let $G$ be a finite group and $\rho:G\rightarrow
U(n)$ be a unitary representation. The $k$-th Newton class
$s_k(\rho)\in H^{2k}(G;\Z)$ of  $\rho$ is
the pull-back of $s_k\in H^{2k}(BU(n);\Z)$ by the
continuous map
$B\rho:BG\rightarrow BU(n)$ induced by $\rho$.
The $k$-th Chern class $c_k(\rho)$ of $\rho$ is defined
similarly. Newton classes of unitary representations
satisfy the following properties:
\begin{enumerate}
\item[(a)] If $\rho_1$ and $\rho_2$ are unitary representations of $G$,
then, for each $k\geq 1$,
\[
s_k(\rho_1\oplus\rho_2)=s_k(\rho_1)+s_k(\rho_2).
\]
\item[(b)] If $\rho:G\rightarrow U(1)$ is a linear character, then,
for each $k\geq 1$ and $l\geq 0$,
\[
s_k(\rho)=s_1(\rho)^k=c_1(\rho)^k,\
s_1(\rho^{\otimes l})=ls_1(\rho).
\]
\end{enumerate}

\section{Fixed point data}\label{sec-fpd}
Let $\gamma\in\Gamma_g$ be an element of order $n$
and $\lgr$ the subgroup
generated by $\gamma$.
By a classical result of Nielsen, $\lgr$ is realized as a group of
automorphisms of a suitable compact Riemann surface $M$ of genus $g$.
Let $h(\gamma)$ be the genus of the quotient surface $M/\lgr$,
and $\lgr_x$ be the isotropy subgroup at $x\in M$.
Set $\calS=\{x\in M\ |\ \lgr_x\not=1\}$ and 
choose a set $\{x_i\}_{1\leq i\leq q}\subset\calS$ of
representatives of $\lgr$-orbits of elements in $\calS$.
Let $\alpha_i$ be the order of the isotropy subgroup at $x_i$
and
$\beta_i$  an integer such that $\gamma^{\beta_i n/\alpha_i}$
acts on $T_{x_i}M$ by
$z\mapsto\exp(2\pi\sqrt{-1}/\alpha_i)z$ with respect to a suitable
local coordinate $z$ at $x_i$.
The number $\beta_i$ is well-defined modulo $\alpha_i$ and is prime to
$\alpha_i$.
The {\em fixed point data} of $\gamma$ is the collection
\[
\langle g,n\ |\ \beta_1/\alpha_1,\cdots,\beta_q/\alpha_q\rangle.
\]
The rational numbers $\beta_1/\alpha_1,\cdots,\beta_q/\alpha_q$,
consider as elements in $\Q/\Z$, are unique up to order.
According to Nielsen \cite{nielsen} (see also \cite{symonds}),
the fixed point data is
independent of various choices made, and hence is
well defined for $\gamma\in\Gamma_g$.
The following proposition will be used in the proof of Theorem
\ref{thm-main}.
\begin{pro}\label{pro-RH}
The fixed point data
$\langle g,n\ |\ \beta_1/\alpha_1,\cdots,\beta_q/\alpha_q\rangle$
of $\gamma\in\Gamma_g$ satisfies
\begin{equation}\label{eq:fpd1}
\sum_{i=1}^q\frac{\beta_i}{\alpha_i}\in\Z
\end{equation}
and the Riemann-Hurwitz equation
\[
2g-2=n(2\cdot h(\gamma)-2)
+n\sum_{i=1}^q\left(1-\frac{1}{\alpha_i}\right).
\]
\end{pro}
\begin{proof} See \cite{yagita}. \end{proof}

As defined in the introduction,
let $e_k(\gamma),s_k(\gamma)\in H^{2k}(\langle\gamma\rangle;\Z)$
be the restriction of $e_k,s_k\in H^{2k}(\Gamma_g;\Z)$ to
$\lgr$.
We wish to express $e_k(\gamma),s_k(\gamma)$ in terms of
the fixed point data of $\gamma$.
Let $\rho_\gamma:\langle\gamma\rangle\rightarrow U(1)$ 
be the linear character defined by
$\gamma\mapsto\exp (2\pi i/n)$ and
$u_{\gamma}:=c_1(\rho_{\gamma})\in
H^2(\langle\gamma\rangle;\Z)$ the first Chern class of
$\rho_\gamma$.
As is known, we have
\[
H^*(\langle\gamma\rangle,\Z)\cong\Z[u_{\gamma}]/(nu_{\gamma}).
\]
Then $e_k(\gamma)$ is expressed in terms of
the fixed point data of $\gamma$, 
as stated below in Proposition \ref{pro:MM-cyclic}.
\begin{pro}\label{pro:MM-cyclic}
Let $\gamma\in\Gamma_g$ be an element of order $n$ having the
fixed point data $\langle g,n\ |\
\beta_1/\alpha_1,\cdots,\beta_q/\alpha_q\rangle$.
Then
\[
e_k(\gamma)=\sum_{i=1}^q
\frac{n}{\alpha_i}(\beta_i^*)^k
u_{\gamma}^k\in H^{2k}(\langle\gamma\rangle;\Z),
\]
where $\beta_i^*$ is an integer satisfying $\beta_i^*\beta_i\equiv
1\pmod{\alpha_i}$.
\end{pro}
\begin{proof}
See \cite{AKU,MR1664755}.
\end{proof}

Let $H^0(M,\Omega_M^1)$ be the space of holomorphic
1-forms of $M$. The natural action of $\lgr$ on
$H^0(M,\Omega_M^1)$ yields a linear representation
$\omega:\lgr\rightarrow U(g)$.
and the induced map
$B\omega:B\lgr\rightarrow BU(g)$ of classifying spaces
can be identified with the composition
$B\lgr\rightarrow B\Gamma_g\rightarrow BU(g)$,
where the former is the map induced by the inclusion
$\lgr\hookrightarrow\Gamma_g$ and the latter is
the map $\eta$ mentioned in Section \ref{sec-def}
(see Remark \ref{rem:holo} in Section \ref{sec-def}).
Hence $s_k(\gamma) \in H^{2k}(\lgr;\Z)$ coincides
with the $k$-th Newton class
$s_k(\omega)$ of the representation $\omega$.
Since $\lgr$ is cyclic,
the representation $\omega$ decomposes into
the sum of linear characters as $\omega=\sum_{j=0}^{n-1}
n_j\rho_{\gamma}^j$.
The multiplicity $n_j$ can be expressed
in terms of the fixed point data,
as in the proposition
below. For each rational number $r\in\Q$, define
$\delta(r)\in\Z$ by
\[
\delta(r)=\begin{cases} 1 & r\in\Z \\
0 & r\not\in\Z. \end{cases}
\]
\begin{pro}\label{pro-multi}
Let $\langle g,n\ |\
\beta_1/\alpha_1,\cdots,\beta_q/\alpha_q\rangle$
be the fixed point data of $\gamma^{-1}\in\Gamma_g$.
Then $n_0=h(\gamma)$ and
\begin{equation}\label{mult-fp}
n_j=h(\gamma)-1+q-\sum_{i=1}^q\left\{
\frac{j\beta_i}{\alpha_i}
-\left[\frac{j\beta_i}{\alpha_i}\right]
+\delta\left(\frac{j\beta_i}{\alpha_i}\right)\right\}
\end{equation}
for $1\leq j<n$.
\end{pro}
To prove the proposition, we follow the argument of
\cite{farkas-kra} first.
Let $\pi:M\rightarrow M/\lgr$ be the natural projection.
For each $1\leq l<n$ with $l|n$, we set
\[
X_l:=\{P\in M\ |\ \gamma^lP=P,\gamma^kP\not=P\ (0<k\leq l-1)\}.
\]
For each non-empty $X_l$, we set
\[
\pi(X_l)=\{Q_{ml}\in M/\lgr\ |\ 1\leq m\leq y_l\},
\]
where $y_l$ is the cardinality of $X_l$
($y_l=0$ when $X_l=\emptyset$).
Choose a local coordinate $z$ for each point in
$\pi^{-1}(Q_{ml})$ such that $\gamma^{-1}$ is given by
% \[
$\gamma^{-1}:z\mapsto\eta_{ml}z$,
% \]
where $\eta_{ml}$ is a primitive $n/l$-root of unity.
For each $j$ and $Q_{ml}\in\pi(X_l)$,
choose the unique natural number
$\lambda_{mlj}$ such that
\begin{equation}\label{cond-lambda}
1\leq\lambda_{mlj}\leq n/l,\
\eta_{ml}^{\lambda_{mlj}}=\epsilon^{jl}.
\end{equation}
Then  (2.5.4) of \cite[Chapter V]{farkas-kra} asserts that
$n_0=h(\gamma)$ and
\begin{equation}\label{eq-FK}
n_j=h(\gamma)-1+\sum_{l|n}y_l-\frac{1}{n}\sum_{l|n}\sum_{m=1}^{y_l}
l\lambda_{mlj}
\end{equation}
for $1\leq j<n$.
\begin{proof}[Proof of Proposition \ref{pro-multi}]
Let $x_i\in M$ be the point corresponding to
$\beta_i/\alpha_i$ $(1\leq i\leq q)$
and suppose that $x_i\in\pi^{-1}(Q_{ml})$.
Then $l=n/\alpha_i$, and the condition \eqref{cond-lambda}
is equivalent to the condition
%\begin{equation}\label{lambda}
\[
1\leq\lambda_{mlj}\leq \alpha_i,\
\lambda_{mlj}\equiv j\beta_i\pmod{\alpha_i},
\]
%\end{equation}
%Then \eqref{lambda} leads to
which leads to
\begin{equation}\label{lambda2}
\lambda_{mlj}=j\beta_i-\alpha_i
\left[\frac{j\beta_i}{\alpha_i}\right]
+\alpha_i\cdot\delta
\left(\frac{j\beta_i}{\alpha_i}\right).
\end{equation}
Substitute \eqref{lambda2} to \eqref{eq-FK}, we have
\begin{align*}
n_j=& h(\gamma)-1+q-\frac{1}{n}\sum_{l|n}\sum_{m=1}^{y_l}l\lambda_{mlj} \\
=& h(\gamma)-1+q-\frac{1}{n}\sum_{l|n}\sum_{
\substack{
1\leq i\leq q \\
\alpha_i=n/l}}
\frac{n}{\alpha_i}
\left\{
j\beta_i-\alpha_i
\left[\frac{j\beta_i}{\alpha_i}\right]
+\alpha_i\cdot\delta
\left(\frac{j\beta_i}{\alpha_i}\right)
\right\} \\
=& h(\gamma)-1+q-\sum_{i=1}^q \left\{
\frac{j\beta_i}{\alpha_i}-\left[\frac{j\beta_i}{\alpha_i}\right]
+\delta\left(\frac{j\beta_i}{\alpha_i}\right)\right\}
\end{align*}
as desired.
\end{proof}
Now we can express $s_k(\gamma)$ in terms of the fixed point data
of $\gamma$, as stated below in Proposition \ref{pro-newton}.
\begin{pro}\label{pro-newton}
Let $\gamma\in\Gamma_g$ be an element of order $n$
having the fixed point data
$\langle g,n\ |\
\beta_1/\alpha_1,\cdots,\beta_q/\alpha_q\rangle$.
For each $j$ $(1\leq j<n)$, set
\[
n_j=h(\gamma)-1+q-\sum_{i=1}^q\left\{
\frac{j\beta_i}{\alpha_i}
-\left[\frac{j\beta_i}{\alpha_i}\right]
+\delta\left(\frac{j\beta_i}{\alpha_i}\right)\right\}
\]
as in \eqref{mult-fp}. Then
\[
s_k(\gamma)=(-1)^k\sum_{j=1}^{n-1}n_j j^ku_{\gamma}^k
\in H^{2n}(\lgr;\Z).
\]
\end{pro}
\begin{proof}
By Proposition \ref{pro-multi},
the representation $\omega:\lgr\rightarrow U(g)$ decomposes
into the sum of linear characters as
$\omega=\sum_{j=0}^{n-1}n_j\rho_{\gamma^{-1}}^j$
with $n_0=h(\gamma)$.
In view of Section \ref{sec-def} (a) and (b), we see that
\[
s_1(\rho_{\gamma^{-1}}^j)=js_1(\rho_{\gamma^{-1}})
=-js_1(\rho_{\gamma})=-ju_{\gamma}
\]
and
%and $s_1(\rho_{\gamma^{-1}}^j)=js_1(\rho_{\gamma^{-1}})$,
\[
s_k(\gamma)=s_k(\omega)
=\sum_{j=0}^{n-1}n_j s_k(\rho_{\gamma^{-1}}^j)
=\sum_{j=1}^{n-1}n_j (-ju_{\gamma})^k
=(-1)^k\sum_{j=1}^{n-1}n_j j^ku_{\gamma}^k
\]
as desired.
\end{proof}

\section{Proof of Theorem \ref{thm-main}}\label{sec-proof}
The proof of Theorem \ref{thm-main} relies on
a number of arithmetic properties of Bernoulli numbers.
We recall those properties first.
Define integers $N_{2k},D_{2k},N_{2k}',D_{2k}'$ by
\[
B_{2k}=\frac{N_{2k}}{D_{2k}},\
\frac{B_{2k}}{2k}=\frac{N_{2k}'}{D_{2k}'},\
(N_{2k},D_{2k})=(N_{2k}',D_{2k}')=1,\
D_{2k},D_{2k}'>0
\]
($N_{2k}',D_{2k}'$ are already introduced in the introduction).
Note that $D_{2k}$ divides $D_{2k}'$.
\begin{lem}\label{lem-divisibility}
Let $p$ be a prime number. Then
\begin{enumerate}
\item $p$ divides $D_{2k}$ if and only if $p-1|2k$,
% \item In particular, $6$ always divides $D_{2k}$.
\item $p^2$ does not divide $D_{2k}$ for all $k\geq 1$.
\item $6$ divides $D_{2k}$ and $4$ divides $D_{2k}'$ for
all $k\geq 1$.
\end{enumerate}
\end{lem}
\begin{proof}
Statements (1) and (2) are direct consequences of
\cite[p.233 Corollary 3]{IR}.
The statement (3) follows from (1) and (2).
\end{proof}
For positive integers $m,n>0$, define
\[
\funcS_m(n):=\sum_{k=1}^{n-1}k^m=1^m+2^m+\cdots+(n-1)^m.
\]
\begin{lem}\label{cong-sum}
Let $m,n$ be positive integers. Then
\[
2\funcS_{2m-1}(n)\equiv 0,\
D_{2m}\funcS_{2m}(n)\equiv 0\pmod{n}. \]
\end{lem}
\begin{proof}
The former follows from
\[
2\sum_{k=1}^{n-1}k^{2m-1}\equiv
\sum_{k=1}^{n-1}k^{2m-1}+\sum_{k=1}^{n-1}(-k)^{2m-1}=0\pmod{n}.
\]
The latter is a direct consequence of \cite[Proposition 15.2.2]{IR}.
\end{proof}
\begin{lem}\label{lem-summation-mod-p}
Let $p$ be a prime number and $b,l$ positive integers. Then
\begin{enumerate}
\item If $l$ is odd then $2\funcS_{l}(p^b)\equiv 0\pmod{p^b}$.
\item If $l$ is even then $p\funcS_{l}(p^b)\equiv 0\pmod{p^b}$.
\end{enumerate}
\end{lem}
\begin{proof}
The first statement is a special case of Lemma \ref{cong-sum}.
The second statement follows from Lemma \ref{lem-divisibility} (2) and
Lemma \ref{cong-sum}.
\end{proof}
Theorem \ref{thm-voronoi} stated in the below,
which is a generalization of the classical Voronoi's congruence
\cite[Proposition 15.2.3]{IR},
is the key ingredient for the proof of
Theorem \ref{thm-main}. The proof of Theorem \ref{thm-voronoi}
is postponed to the next section.
\begin{thm}\label{thm-voronoi}
Let $p$ be a prime number and $c$ an integer prime to $p$.
Then
\begin{align}
N_{2k}'\cdot p^b(c^{2k}-1)\equiv
D_{2k}'\cdot c^{2k-1}\sum_{s=0}^{p^{a+b}-1}
s^{2k-1}\left[\frac{sc}{p^a}\right]
\pmod{p^{a+b}}
\end{align}
for all $a,b\geq 0$ and $k\geq 1$.
\end{thm}
Now we prove Theorem \ref{thm-main}.
Since $e_k(\gamma)$ and $s_k(\gamma)$ are natural with respect to
the restriction to subgroups,
it suffices to prove the case when $\lgr$ is a $p$-group.
Let $p$ be a prime and $\gamma\in\Gamma_g$ an element
of order $p^n$ having the fixed point data
\[
\langle g,p^n\ |\ \beta_1/p^{\alpha_1},
\cdots,\beta_q/p^{\alpha_q}\rangle
\]
(we write $p^{\alpha_i}$ instead of $\alpha_i$).
Choose an integer
$\beta_i^*$ satisfying $\beta_i^*\beta_i\equiv
1\pmod{p^n}$ for each $1\leq i\leq q$,
and set
\[
n_j=h(\gamma)-1+q-\sum_{i=1}^q \left\{
\frac{j\beta_i}{p^{\alpha_i}}-\left[\frac{j\beta_i}{p^{\alpha_i}}\right]
+\delta\left(\frac{j\beta_i}{p^{\alpha_i}}\right)\right\}
\]
for each $1\leq j\leq p^n-1$.
According Proposition \ref{pro:MM-cyclic} and
Proposition \ref{pro-newton}, we have
\[
e_{2k-1}(\gamma)=\sum_{i=1}^qp^{n-\alpha_i}
(\beta_i^*)^{2k-1}u_{\gamma}^{2k-1},\ 
s_{2k-1}(\gamma)=-\sum_{j=1}^{p^n-1}n_jj^{2k-1}u_{\gamma}^{2k-1}.
\]
Hence it suffices to prove the congruence
\begin{equation}\label{want}
N_{2k}'\sum_{i=1}^qp^{n-\alpha_i}(\beta_i^*)^{2k-1}
\equiv -D_{2k}'\sum_{j=1}^{p^n-1}n_jj^{2k-1}
\pmod{p^n}.
\end{equation}
Since $2\sum_{j=1}^{p^n-1}j^{2k-1}\equiv0\pmod{p^n}$ by
Lemma \ref{cong-sum} and $D_{2k}'$ is even
by Lemma \ref{lem-divisibility},
the right hand side of \eqref{want} satisfies
\begin{align*}
(\mathrm{RHS})=& \ -D_{2k}'\sum_{j=1}^{p^n-1}j^{2k-1}\left[
h-1+q-\sum_{i=1}^q \left\{
\frac{j\beta_i}{p^{\alpha_i}}-\left[\frac{j\beta_i}{p^{\alpha_i}}\right]
+\delta\left(\frac{j\beta_i}{p^{\alpha_i}}\right)\right\}
\right] \\
\equiv & \ D_{2k}'\sum_{j=1}^{p^n-1}j^{2k-1}
\sum_{i=1}^q \left\{
\frac{j\beta_i}{p^{\alpha_i}}-\left[\frac{j\beta_i}{p^{\alpha_i}}\right]
+\delta\left(\frac{j\beta_i}{p^{\alpha_i}}\right)\right\} 
\pmod{p^n} \\
=&\ D_{2k}'\left\{
\sum_{i=1}^q\frac{\beta_i}{p^{\alpha_i}}\sum_{j=1}^{p^n-1}j^{2k}
-\sum_{i=1}^q\sum_{j=1}^{p^n-1}j^{2k-1}
    \left[\frac{j\beta_i}{p^{\alpha_i}}\right]
+\sum_{i=1}^q\sum_{j=1}^{p^n-1}j^{2k-1}
    \delta\left(\frac{j\beta_i}{p^{\alpha_i}}\right)
\right\}.
\end{align*}
Since $\sum_{i=1}^q\beta_i/p^{\alpha_i}\in\Z$ by Proposition
\ref{pro-RH} \eqref{eq:fpd1}, and 
$D_{2k}'\sum_{j=1}^{p^n-1}j^{2k}\equiv 0\pmod {p^n}$
by Lemma \ref{cong-sum}, we see that
\[
D_{2k}'\sum_{i=1}^q
\frac{\beta_i}{p^{\alpha_i}}\sum_{j=1}^{p^n-1}j^{2k}
\equiv 0\pmod{p^n}.
\]
On the other hand, we have
\[
D_{2k}'\sum_{j=1}^{p^n-1}j^{2k-1}\delta
\left(\frac{j\beta_i}{p^{\alpha_i}}\right)
=D_{2k}'\sum_{\substack{
   1\leq j\leq p^n-1 \\ p^{\alpha_i}|j}}j^{2k-1}
\equiv 0\pmod {p^n},
\]
where the first equality follows from the fact $\beta_i$ is prime to $p$,
while the second congruence follows from the fact $D_{2k}'$ is even
as in the proof of Lemma \ref{cong-sum}.
Hence the right hand side of \eqref{want} satisfies
\begin{equation}\label{newton2}
(\mathrm{RHS})=
-D_{2k}'\sum_{j=1}^{p^n-1}n_jj^{2k-1}\equiv
-D_{2k}'\sum_{i=1}^q\sum_{j=1}^{p^n-1}j^{2k-1}
\left[\frac{j\beta_i}{p^{\alpha_i}}\right]\pmod {p^n}.
\end{equation}
Now Theorem \ref{thm-voronoi} implies
\[
N_{2k}'p^{n-\alpha_i}(\beta_i^{2k}-1)\equiv
D_{2k}'\beta_i^{2k-1}\sum_{j=1}^{p^n-1}j^{2k-1}
\left[\frac{j\beta_i}{p^{\alpha_i}}\right]\pmod {p^n}.
\]
Multiplying $-(\beta_i^*)^{2k-1}$ to the both sides of
the last congruence, we have
\begin{equation}\label{eq:pre-sum}
N_{2k}'p^{n-\alpha_i}((\beta_i^*)^{2k-1}-\beta_i)\equiv
-D_{2k}'\sum_{j=1}^{p^n-1}j^{2k-1}
\left[\frac{j\beta_i}{p^{\alpha_i}}\right]\pmod {p^n}.
\end{equation}
Since
\[
\sum_{i=1}^q p^{n-\alpha_i}\beta_i
=p^n\sum_{i=1}^q \frac{\beta_i}{p^{\alpha_i}}\equiv 0\pmod {p^n}
\]
by Proposition \ref{pro-RH} \eqref{eq:fpd1},
the summation of \eqref{eq:pre-sum} for $1\leq i\leq q$
provides a congruence
\[
N_{2k}'\sum_{i=1}^q p^{n-\alpha_i}(\beta_i^*)^{2k-1}
\equiv 
-D_{2k}'\sum_{i=1}^q\sum_{j=1}^{p^n-1}j^{2k-1}
\left[\frac{j\beta_i}{p^{\alpha_i}}\right]
%-D_{2k}'\sum_{j=1}^{p^n-1}n_jj^{2k-1}
\pmod {p^n}.
\]
In view of \eqref{newton2}, this is the desired
congruence \eqref{want}, and hence
verifying Theorem \ref{thm-main}.

\section{A generalization of Voronoi's congruence
}\label{sec-Voronoi}
We prove Theorem \ref{thm-voronoi} in this section.
The proof of Theorem \ref{thm-voronoi} is based on
a generalization of Voronoi's congruences due to Porubsk\'y
stated below.
\begin{thm}[Porubsk\'y \cite{porubsky,porubsky2}]\label{thm-por}
Given a positive integer $N$ and
%a rational number $c$ prime to $N$
%$($i.e. the denominator and the numerator of $c$ are prime to $N)$,
an integer $c$ prime to $N$,
then
\begin{align}\label{eq-por1}
& (c^{2k}-1)\frac{B_{2k}}{2k}+\frac{c^{2k}-c^{2k-1}}{2}\cdot\frac{2k-1}{2}
B_{2k-2}N^2 \notag \\
\equiv & c^{2k-1}\sum_{x=1}^{N-1}x^{2k-1}
\left[\frac{cx}{N}\right]\pmod{N}
\end{align}
for all $k\geq 1$. In particular, if $N$ is odd or $N(c-1)\equiv 0
\pmod{8}$, then we have
\begin{equation}\label{eq-por2}
(c^{2k}-1)\frac{B_{2k}}{2k}\equiv
c^{2k-1}\sum_{x=1}^{N-1}x^{2k-1}
\left[\frac{cx}{N}\right]\pmod{N}
\end{equation}
for all $k\geq 1$.
\end{thm}
Here congruences in Theorem \ref{thm-por}
%\eqref{eq-por1} and \eqref{eq-por2}
are those in $\Z_{(N)}=\{m/n\in\Q\ |\ (n,N)=1\}$.
Namely, for $r,s\in\Z_{(N)}$, we write $r\equiv s\pmod{N}$
if $r-s=m/n, (n,N)=1$, and $N|m$.
\begin{rem}
It was claimed in \cite{porubsky} that the congruence
\eqref{eq-por1} holds for any rational number $c$ prime to $N$.
However, as was corrected in \cite{porubsky2}, the congruence
is valid only when $c$ is an integer.
\end{rem}
If $a=0$ then Theorem \ref{thm-voronoi} is a direct consequence
of Lemma \ref{cong-sum}. Henceforth, we assume $a>0$.
We first prove Theorem \ref{thm-voronoi} for the case $b=0$.
\begin{pro}\label{pro-coprime}
Let $p$ be a prime number and $c$ an integer prime to $p$.
Then
\begin{equation*}
N_{2k}'\cdot (c^{2k}-1)\equiv
D_{2k}'\cdot c^{2k-1}\sum_{x=1}^{p^a-1}x^{2k-1}
\left[\frac{cx}{p^a}\right]\pmod{p^a}
\end{equation*}
for all $a,k\geq 1$.
\end{pro}
\begin{proof}
The proposition is a direct consequence of Theorem \ref{thm-por}
\eqref{eq-por2}
unless $p^a=2$. Suppose $p=2$ and $a=1$. Then we have
\begin{align*}
\frac{c^{2k}-c^{2k-1}}{2}\cdot\frac{2k-1}{2} B_{2k-2}\cdot 2^2
=(c^{2k}-c^{2k-1})(2k-1)B_{2k-2}\in\Z_{(2)},
\end{align*}
because $c^{2k}-c^{2k-1}$ is even and $2B_{2k-2}\in\Z_{(2)}$
by Lemma \ref{lem-divisibility} (2).
Since $D_{2k}'$ is even, we have
\[
D_{2k}'\cdot\frac{c^{2k}-c^{2k-1}}{2}\cdot\frac{2k-1}{2} B_{2k-2}\cdot 2^2
%D_{2k}'(c^{2k}-c^{2k-1})(2k-1) B_{2k-2}
\equiv 0\pmod{2}.
\]
Now the result follows from Theorem \ref{thm-por} \eqref{eq-por1}.
\end{proof}\noindent
Now we prove Theorem \ref{thm-voronoi}.
According to Proposition \ref{pro-coprime}, we have
\[
N_{2k}'\cdot p^b(c^{2k}-1)
\equiv  D_{2k}'\cdot p^b c^{2k-1}\sum_{x=1}^{p^a-1}x^{2k-1}
\left[\frac{cx}{p^a}\right]\pmod{p^{a+b}}.
\]
Hence it suffices to prove
\begin{equation}\label{eq-intermediate-cong}
D_{2k}'\cdot
p^b \sum_{s=0}^{p^a-1}s^{2k-1}
\left[\frac{sc}{p^a}\right]
\equiv
D_{2k}' 
\sum_{s=0}^{p^{a+b}-1}s^{2k-1}\left[\frac{sc}{p^a}\right]
\pmod{p^{a+b}}.
\end{equation}
Now since
\begin{align*}
& \sum_{s=0}^{p^{a+b}-1}s^{2k-1}\left[\frac{sc}{p^a}\right]
=\sum_{j=0}^{p^b-1}\sum_{s=0}^{p^a-1}(s+jp^a)^{2k-1}
\left[\frac{(s+jp^a)c}{p^a}\right], \\
& \left[\frac{(s+jp^a)c}{p^a}\right]=jc+\left[\frac{sc}{p^a}\right],
\end{align*}
we see that
\[
\sum_{s=0}^{p^{a+b}-1}s^{2k-1}\left[\frac{sc}{p^a}\right]
=\sum_{j=0}^{p^b-1}\sum_{s=0}^{p^a-1}\left(
cj(s+jp^a)^{2k-1}+(s+jp^a)^{2k-1}
\left[\frac{sc}{p^a}\right]\right).
\]
\begin{lem}\label{lem-vanish1}
Under the assumption of Theorem \ref{thm-voronoi}, we have
\[
D_{2k}'
\sum_{j=0}^{p^b-1}\sum_{s=0}^{p^a-1}j(s+jp^a)^{2k-1}\equiv 0
\pmod{p^{a+b}}.
\]
\end{lem}
\begin{proof}
Observe first
\begin{align*}
\sum_{j=0}^{p^b-1}\sum_{s=0}^{p^a-1}j(s+jp^a)^{2k-1}
=& \sum_{j=0}^{p^b-1}\sum_{s=0}^{p^a-1}
\sum_{l=0}^{2k-1}j\binom{2k-1}{l}s^{2k-1-l}(jp^a)^l \\
=& \sum_{j=0}^{p^b-1}\sum_{l=0}^{2k-1}\binom{2k-1}{l}
j^{l+1}p^{al}\cdot\funcS_{2k-1-l}(p^a) \\
=& \sum_{l=0}^{2k-1}\binom{2k-1}{l}p^{al}\cdot
\funcS_{2k-1-l}(p^a)\funcS_{l+1}(p^b).
\end{align*}
In view of Lemma \ref{lem-summation-mod-p} and the fact
that $D_{2k}$ is even, we see that
\begin{equation}\label{eq-interm}
D_{2k}\cdot p^{al}\cdot \funcS_{2k-1-l}(p^a)\funcS_{l+1}(p^b)
\equiv 0\pmod{p^{a+b}}
\end{equation}
provided $al\geq a+1$ (i.e. $l\geq 2$).
Now $D_2\cdot\funcS_2(p^b)=6\funcS_2(p^b)\equiv 0\pmod{p^b}$.
Since $D_{2k}$ is divisible by 6 by Lemma \ref{lem-divisibility},
we see that
\begin{equation}
D_{2k}\cdot p^{a}\cdot \funcS_{2k-2}(p^a)\funcS_{2}(p^b)
\equiv 0\pmod{p^{a+b}},
\end{equation}
which shows that the congruence \eqref{eq-interm} is
valid for $l=1$.
For $l=0$, recall that $D_{2k}'$ is divisible by 4
by Lemma \ref{lem-divisibility} so that
$D_{2k}'=2\gamma D_{2k}$ for some positive integer $\gamma$.
Hence
\[
D_{2k}'\cdot
\funcS_{2k-1}(p^a)\funcS_{1}(p^b)
=\gamma D_{2k}\cdot
\funcS_{2k-1}(p^a)\cdot p^b(p^b-1)\equiv 0\pmod{p^{a+b}}
\]
as $D_{2k}$ is even and $2\funcS_{2k-1}(p^a)\equiv 0\pmod{p^a}$.
Putting all these together, we obtain
\[
D_{2k}'\sum_{l=0}^{2k-1}\binom{2k-1}{l}p^{al}\cdot
\funcS_{2k-1-l}(p^a)\funcS_{l+1}(p^b)\equiv 0\pmod{p^{a+b}},
\]
which prove the lemma.
\end{proof}
\begin{lem}\label{lem-vanish2}
Under the assumption of Theorem \ref{thm-voronoi}, we have
\begin{align*}
D_{2k}'\sum_{j=0}^{p^b-1}\sum_{s=0}^{p^a-1}
(s+jp^a)^{2k-1}\left[\frac{sc}{p^a}\right]
\equiv
D_{2k}'\cdot p^b \sum_{s=0}^{p^a-1}s^{2k-1}
\left[\frac{sc}{p^a}\right]\pmod{p^{a+b}}.
\end{align*}
\end{lem}
\begin{proof}
As in the proof of the previous lemma, we have
\[
\sum_{j=0}^{p^b-1}\sum_{s=0}^{p^a-1}(s+jp^a)^{2k-1}
\left[\frac{sc}{p^a}\right]
=\sum_{l=0}^{2k-1}\binom{2k-1}{l}p^{al}\funcS_l(p^b)
\sum_{s=0}^{p^a-1}s^{2k-1-l}\left[\frac{sc}{p^a}\right].
\]
In view of Lemma \ref{lem-summation-mod-p}, we have
$D_{2k}\cdot p^{al}\funcS_l(p^b)\equiv 0\pmod{p^{a+b}}$
if $l$ is odd or if $l$ is even and $al\geq a+1$ (i.e. $l\not=0$).
Thus
\begin{align*}
& \ D_{2k}\sum_{l=0}^{2k-1}\binom{2k-1}{l}p^{al}\funcS_l(p^b)
\sum_{s=0}^{p^a-1}s^{2k-1-l}\left[\frac{sc}{p^a}\right] \\
\equiv & \ D_{2k}\funcS_0(p^b)\sum_{s=0}^{p^a-1}s^{2k-1}
\left[\frac{sc}{p^a}\right]
=D_{2k} p^b \sum_{s=0}^{p^a-1}s^{2k-1}
\left[\frac{sc}{p^a}\right]\pmod{p^{a+b}}.
\end{align*}
The lemma follows immediately.
\end{proof}\noindent
Combining Lemma \ref{lem-vanish1} and Lemma \ref{lem-vanish2},
the congruence \eqref{eq-intermediate-cong} follows,
and hence verifying Theorem \ref{thm-voronoi}.

\begin{ack}
The authors thank Prof. \v Stefan Porubsk\'y for infoming
the paper \cite{porubsky2}.
The first author is partially
supported by the Grant-in-Aid for Scientific Research (C)
(No. 17560054) from the Japan Society for Promotion of Sciences.
%and Grant-in-Aid for formation of COE
%in Hokkaido University.
The second author is partially
supported by the Grant-in-Aid for Scientific Research (A)
(No. 18204002) from the Japan Society for Promotion of Sciences.
\end{ack}

\bibliographystyle{amsplain}

\providecommand{\bysame}{\leavevmode\hbox to3em{\hrulefill}\thinspace}

\end{document}